\documentclass[10pt]{amsart}

\usepackage{amssymb}
\usepackage{bm}
\usepackage{graphicx}
\usepackage[centertags]{amsmath}
\usepackage{amsfonts}
\usepackage{amsthm}

\newtheorem{theorem}{Theorem}
\newtheorem*{mtheo}{Theorem}

\theoremstyle{definition}

\newtheorem{claim}{Claim}

\begin{document}
\title{Canonical Borel equivalence relations on $ \mathbb{R}^n$}

\author{VLITAS Dimitris}

\footnotetext{The research leading to these results has received funding from the [European Community's] Seventh Framework Programme [FP7/2007-2013] under grant agreement n 238381}
 %------------------------Abstract-------------------------------%
\begin{abstract}  In this paper we complete the attempt of H. Lefmann to show that Borel equivalence relations on the $n$-element subsets of $2^{\omega}$, that respect an order type, have a finite Ramsey basis.
\end{abstract}

\maketitle

\section{Introduction} The purpose of this paper is to prove a canonical Ramsey theorem for the finite powers of the Cantor space $2^{\omega}$ completing a previous attempt of H.Lefmann \cite{Le}. We shall, of course, use some of the ideas from \cite{Le} but we shall also correct an essential omission of that paper.\\
It is well known how Ramsey Theorem, \cite{Ra} generalizes to its canonical form the Erd\H{o}s--Rado Theorem, \cite{Er--Ra}.  Namely, the classical Ramsey theorem considers finite colorings of $n$-element subsets of $\omega$ and the Erd\H{o}s--Rado theorem considers countable colorings of the same set.

 \begin{mtheo} (Ramsey)\cite{Ra} For every positive integer $n$ and every finite coloring of the family $[\omega]^n$ of all $n$-element subsets of $\omega$, there is an infinite subset $X$ of $\omega$ such that the set $[X]^n$ of all $n$-element subsets of $X$ is monochromatic. \end{mtheo}
 \begin{mtheo}(Erd\H{o}s--Rado)\cite{Er--Ra} Given an arbitrary set $R$, a positive integer $n$ and a mapping $c:[\omega]^n\to R$, there are an infinite subset $X\subseteq \omega$ and a finite set $I\subseteq n$ such that for any $A,B\in [X]^{n}$ one has $c(A)=c(B)$ if and only if $$A:I=B:I$$
\end{mtheo}

Here for $A=(x_0,\dots, x_{n-1})\in [\omega]^n$, $I\subseteq n$ $A:I=\{ x_i:i\in I\}$.

In this paper we consider uncountable versions of these results. In particular we consider the Cantor space $2^{\omega}$ with its lexicographic ordering induced by $0<1$. By $[2^\omega]^2$ we mean the set of all pairs ordered increasingly with respect to the lexicographic order.
In \cite{Sie} it is shown that there exists a coloring $c:[2^{\omega}]^2\to \{0,1\}$ such that for any uncountable subset $Y\subseteq 2^{\omega}$ the restriction of $c$ to $[Y]^2$ is not constant.
The definition of that coloring uses a well-ordering of the continum, so it is not definable. Similarly in \cite{Ga--Sel} is shown that there is a partition of $[2^{\omega}]^2$ into infinitely many pieces, so that for any subset  $X$ of cardinality of the continuum, $[X]^2$ intersects all of the pieces. This coloring also uses a well-ordering of the continum so it is not definable.
 It turns out that if we consider mappings which are Baire measurable a partition result can be obtained.

 We start with the Cantor space $2^\omega$. We consider any $x\in 2^\omega$ as an $\omega$-sequence of $0,1$. Then $2^\omega$ can be endowed with the metric $d$ defined by $$d(x,y)=1/{k+1}\text{, }k=\min\{n:x(n)\neq y(n)\}$$ This metric gives on the Cantor space the usual Tychonoff product topology. In this paper we only consider $n$-element subsets of $2^{\omega}$ ordered lexicographically and by $(x_0,\dots,x_{n-1})_{<_{lex}}$ we express that $x_0<_{lex}x_1<_{lex} \dots <_{lex} x_{n-1}$. In that manner $[2^{\omega}]^n$ the set of all $n$-element subsets of $2^{\omega}$ ordered lexicographically, is a subset of the finite product $(2^{\omega})^n$ and it is a topological space with the subspace topology.

     Let $n\in \omega$ and $(\{ 1,\dots, n-1\}, \preceq)$ be a total order. Then\\ $[2^{\omega}]^n_{\preceq} = \{\, (x_0, \dots, x_{n-1})_{<_{lex}} \in [2^{\omega}]^n : (\forall i\neq j, i,j\in n\setminus \{0\})[d(x_{i-1},x_i)> d(x_{j-1},x_j)$ iff $i\preceq j ] \, \}$. Any two elements $X,Y$ of $[2^\omega]^n_\preceq$ are said to have  the same {\it order type}.
     Consider a Borel equivalence relation $E$ on $[2^\omega]^n_\preceq$. In particular $E$ viewed as a subset of $(2^\omega)^{2n}$ has the property of Baire. Mycielski \cite{Myc} has shown that for any meager subset of $[2^\omega]^n$ there exists a perfect subset $\mathcal{P}$ such that $[\mathcal{P}]^n$ avoids it.  Recall that a subset $\mathcal{P}$ of $2^\omega$ is perfect if it is non empty, closed and has no isolated points. Therefore we can assume that $E$ as a subset of $(2^\omega)^{2n}$ is open. As a consequence we can assume that each equivalence class is open. Define $c:[\mathcal{P}]^n_{\preceq}\to 2^\omega$ by letting $$c(x_0,\dots, x_{n-1})=\min([(x_0,\dots, x_{n-1})]_{E})$$ The map $c$ picks from each equivalence class a minimal element with respect the coordinatewise partial ordering, where on each coordinate we take the lexicographic ordering. Then $c$ is a Borel map representing $E$ on $[\mathcal{P}]^n_{\preceq}$.
 Thus any Borel equivalence relation $E$ on $[2^\omega]^n_\preceq$ can be reduced on a perfect subset $\mathcal{P}$ to a smooth equivalence relation, meaning that there exists a Borel mapping $c$ such that $xEy$ if and only if $c(x)=c(y)$.

  Given any $X= (x_0, \dots, x_{n-1})_{<_{lex}}$ element of ${[2^{\omega}]}^n_\preceq$ let $$D(X)=\{\, x_{i-1}\cap x_{i} : i\in n\setminus \{0\}\, \}$$ where $x_0\cap x_1=\max \{t: t\sqsubseteq x_0$ and $t\sqsubseteq x_1\}$. By $t\sqsubseteq x_0$ we denote that $t$ is an initial segment of $x_0$ and maximality is taken under inclusion. Given $t\sqsubseteq x$, $|t|$ denotes the length of $t$, namely the level on which it lies on $2^{<\omega}$, the binary tree ordered by inclusion. For a $t\in 2^{<\omega}$ by $t^\frown d$, for $d\in \{0,1\}$, we denote the unique extension of $t$ so that $|t^\frown d|=|t|+1$.
  Let  $$\widehat{D(X)}=D(X)\cup \{\, x_i\upharpoonright| x_{j-1}\cap x_{j}|: i\neq j, j\in n\setminus \{0\},i\in n\,\}$$ Then $\widehat{D(X)}$ determines uniquely a rooted subtree of $2^{< \omega}$. Notice that for any $X\in [2^\omega]^n_{\preceq}$, $\widehat{D(X)}$ is the same up to isomorphism. Consider $X=(x_0, \dots, x_{n-1})_{<_{lex}},Y=(y_0, \dots, y_{n-1})_{<_{lex}} \in {[2^{\omega}]}^n_\preceq$ then $g: \widehat{D(X)}\to \widehat{D(Y)}$ defined by: $$g(x_{i-1}\cap x_i)=(y_{i-1}\cap y_i)\text{ and }g(x_i\upharpoonright |x_{j-1}\cap x_j|)=y_i\upharpoonright |y_{j-1}\cap y_j|$$ witness that $\widehat{D(X)}$, $\widehat{D(Y)}$ are isomorphic. Therefore we can pick an element\\ $C^n_{\preceq}= \widehat{D(X_0)}$ as a representative of the class $\{\,\widehat{D(X)}: X\in [2^{\omega}]^n_{\preceq}\,\}$ and for every $X\in {[2^\omega]}^n_\preceq$ an isomorphism $g_X:\widehat{D(X)}\to C^n_\preceq$.	For $I\subseteq \{\, 0, \dots, n-1\, \}$ let $$X:I=\{\, x_i: i\in I\,\}$$ and similarly for $J\subseteq C^n_{\preceq}$ we define the set $$\widehat{D(X)}:J=\{ x\upharpoonright k: g_X(x\upharpoonright k)\in J\}$$ where $C^n_{\preceq}$ is the fixed representative of the class $\{\,\widehat{D(X)}: X\in [2^{\omega}]^n_{\preceq}\,\}$.\\

  We recall the following result of F. Galvin

  \begin{theorem} \cite{Ga} Let $n\leq 3$, $l \in \omega$ and let $\mathcal{T}\subseteq 2^{\omega}$ be a perfect subset. Then for any Baire measurable coloring  $c:[\mathcal{T}]^n \to l$ there exists a perfect subset $\mathcal{P}\subseteq \mathcal{T}$ such that $c\upharpoonright [\mathcal{P}]^n_{\preceq}$ is constant for every total order $(\{1,\dots, n-1\}, \preceq)$.
  \end{theorem}
 
   The general case for any $n\in \omega$ was obtained by A. Blass
   \begin{theorem} \cite{Bla} Let $n,l\in \omega$  and let $\mathcal{T}\subseteq 2^{\omega}$ be a perfect subset. Then for any Baire measurable coloring  $c:[\mathcal{T}]^n \to l$ there exists a perfect subset $\mathcal{P}\subseteq \mathcal{T}$ such that $c\upharpoonright [\mathcal{P}]^n_{\preceq}$ is constant for every total order $(\{1,\dots, n-1\}, \preceq)$.
\end{theorem}
   A subset $\mathcal{P}$ of $2^\omega$ is skew if and only if for every $x_0,x_1,y_0,y_1\in 2^\omega$ with $x_0\neq x_1$ and $y_0\neq y_1$, $d(x_0,x_1)= d(y_0,y_1)$ implies that $x_0\cap x_1=y_0\cap y_1$. It is easily seen (for a proof see  \cite{Bla}) that every perfect subset of $2^\omega$ contains a skew perfect subset.

  \section{Main theorem}

  The main theorem of this section is the following:

   \begin{theorem}
   Let $n\in \omega$, $(\{1,\dots, n-1\},\preceq)$ a total order and $E$ a Borel equivalence relation on $[2^\omega]^n_\preceq$. There exists a skew perfect subset $\mathcal{P}$ of $2^\omega$, a subset $I\subseteq n$ and $J\subseteq C^n_\preceq$ such that for all $X,Y\in [\mathcal{P}]^n_\preceq$ one has $XEY$ if and only if $X:I=Y:I$ and $\widehat{D(X)}:J=\widehat{D(Y)}:J$.
   \end{theorem}

   This theorem tells us that any Borel equivalence relation on $[2^\omega]^n_\preceq$ corresponds to a pair $(I,J)$ such that $I\subseteq n$ and $J\subseteq C^n_{\preceq}$. Notice that there exists only finitely many such a pairs. Therefore we say that Borel equivalence relations on $[2^\omega]^n_\preceq$ have a finite {\it{Ramsey basis}}.

   We prove the following version of Theorem $3$

 \begin{theorem} Let $n\in \omega$ and $(\{ 1,\dots, n-1\}, \preceq)$ be a total order. Further let $\mathcal{T}\subseteq 2^{\omega}$ be a perfect subset and $M$ a metric space. Then for any Baire measurable mapping $c: [\mathcal{T}]^{n}_{\leq}\to M$ there exists a skew perfect subset $\mathcal{P}\subseteq \mathcal{T}$, $I\subseteq \{ 0, \dots, n-1\, \}$ and $J\subseteq C^n_{\preceq}$ such that  for every $X,Y \in [\mathcal{P}]^n_{\preceq}$:\\

 $c(X)=c(Y)$ iff  $X:I=Y:I$ and $\widehat{D(X)}:J=\widehat{D(Y)}:J$.
 \end{theorem}

  \begin{proof}	The proof is done by induction on $n$.\\ First of all notice that $c$ can be assumed to be continuous. For the purpose of Theorem $3$ the metric space $M$ is a Polish space (recall the map $c:[\mathcal{P}]^n_{\preceq}\to 2^\omega$ from the introduction). Theorem $4$ deals with any metric space $M$.  Given now any Baire measurable mapping $c:[\mathcal{R}]^n_\preceq \to M$ where $M$ is a metric space and $\mathcal{R}$ a perfect subset of $2^\omega$, there is always a perfect subset $\mathcal{T}$ of $\mathcal{R}$ such that the restriction of $c$ on $[\mathcal{T}]^n_\preceq$ is continuous. To see this notice that by assuming that $M$ is a separable space the map $c$ can be assumed to be continuous on $[\mathcal{R}]^n_\preceq$ modulo a meager set. In \cite{EFK} is shown that in the case of $[\mathcal{R}]^n_\preceq$ the condition of $M$ being separable can be omitted. Then by Mycielski \cite{Myc}  there is a perfect subset $\mathcal{T}$ so that $[\mathcal{T}]^n_\preceq$ avoids that meager set. Therefore the map $c$ restricted on $[\mathcal{T}]^n_\preceq$ is continuous.

    Let $n=1$ and  $c:[\mathcal{T}]^1\to M$ be a continuous map. This induces another continuous map $c^{\star}:[\mathcal{T}]^2\to \{0,1\}$ defined by $c^{\star}(X,Y)=0$ if $c(X)=c(Y)$ and equal to $1$ otherwise. By Theorem $1$ there exists a skew perfect subset $\mathcal{P}\subseteq \mathcal{T}$ and $l\in 2$ such that $c^{\star}\upharpoonright [\mathcal{P}]^2=\{l\}$. In other words the restriction of $c$ on $\mathcal{P}$ is either constant or one to one. In the first case $I=\emptyset$, $J=\emptyset$ and in the second case $I=n$ and $J=\emptyset$.

 To establish the inductive step we need some special constructions. Let $X\in [\mathcal{T}]^n_{\preceq}$ with $X=(x_0, \dots, x_{n-1})_{<_{lex}}$. Consider the map \\ $f:\{0,\dots, n-1\}\to \{1,-1\}$ defined by: for $i\neq 0, n-1$  \begin{equation*} f(i)=\left\{
                             \begin{array}{rl}
                             -1 & \text{if  } i+1\preceq i,\\
                             1 & \text{if  } i\preceq i+1.
                             \end{array}\right.
                             \end{equation*}
for $i=n-1$ $f(i)=-1$ and for $i=0$ $f(i)=1$.
                           \begin{figure}[htb]
                            \centering \includegraphics [width=.65\textwidth]{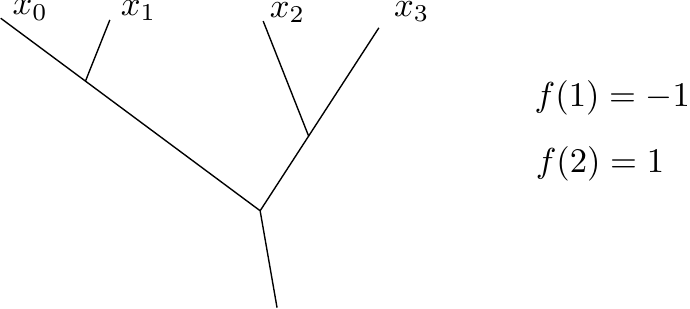}

                             \end{figure}

    In other words $f(i)$, for $i\in n$ and $i\neq0, n-1$, is the unique $j$, with $j\in \{1,-1\}$, such that $x_i\cap x_{i-j}$ is an initial segment of $x_i\cap x_{i+j}$.

    By $X+i\subset [\mathcal{T}]^{n+1}$ we define the set  of $n+1$-tuples resulting from $X$ with the addition of a new element $x^{+}_i$ such that the following conditions hold:\\
    $\bf{1}$ $x_{i-1}<_{lex} x^{+}_i<_{lex} x_i$\\
    $\bf{2}$ $x_i\cap x^{+}_i$ is an initial segment of $x_i\cap x_{i+f(i)}$\\
    $\bf{3}$ $X$ and $(X\setminus \{x_i\})\cup \{x^{+}_i\}$ have the same order-type.\\

    Notice that all elements of $ X+i$ have the same order type $(\{ 1,\dots, n\}, \preceq^{+}_i)$, cause by conditions $\bf{2}$ and $\bf{3}$ one has that $i\in \{1, \dots, n\}$ must have a fixed position in $(\{1,\dots, n\}, \preceq^+_i)$. To see that, below we exhaust all possibilities.\\
    
    Observe that the $i\in ( \{1, \dots, n-1\},\preceq)$ becomes $i+1\in ( \{1, \dots, n\},\preceq^+)$. The $i$ in $ (\{1, \dots, n\},\preceq^+)$ comes from $x_i^+\cap x_{i-1}$ and the $i+1\in ( \{1, \dots, n\},\preceq^+)$ comes from $x_i^+\cap x_{i}$. Keeping that in mind we consider the following two possibilities:
    
    Let $i,j\in (\{1, \dots, n-1\},\preceq)$ and $x_j<_{lex}x_i$. The new element $x_i^+$ also satisfies $x_j<_{lex}x_i^+$.
    If $i\preceq j$, then $i\preceq^+_i j$, for $i,j\in (\{1, \dots, n\},\preceq^+_i)$.
    If $j\preceq i$, then $j\preceq^{+}_i i$ for  $i,j\in (\{1, \dots, n\},\preceq^+_i)$.

     Let $i,j\in (\{1, \dots, n-1\},\preceq)$ and $x_i<_{lex} x_j$. Then $x_i^+ <_{lex} x_j$ as well.
     If $i\preceq j$, then $i\preceq^+_i j+1$, for $i,j+1\in (\{1, \dots, n\},\preceq^+_i)$. In the case of  $j\preceq i $, then $j+1\preceq^{+}_i i$, for $j+1,i\in (\{1,\dots, n\},\preceq^+_i)$. Notice that since $x_i<_{lex} x_j$, $j\in (\{1, \dots, n-1\},\preceq)$ becomes $j+1\in (\{1,\dots, n\},\preceq^+_i)$ (see figures below).\\

    By $X\oplus i\subset [\mathcal{T}]^{n+1}$ we define the set of $n+1$-tuples that result from $X$ with the addition of a new element $x^{\oplus}_i$ such that the following three conditions hold:\\
    $\bf{1}$ $x_{i-1}<_{lex} x^{\oplus}_i<_{lex} x_i$\\
    $\bf{2}$ $x_{i+f(i)} \cap x_i$ is a proper initial segment of $x^{\oplus}_i\cap x_i$.\\
    $\bf{3}$ $X$ and $(X\setminus \{x_i\})\cup \{x^{\oplus}_i\}$ have the same order type.\\

    Notice that all elements of $X\oplus i$, are not necessarily of the same order type, because of condition $\bf{2}$ above $i \in \{ 1,\dots, n\}$ cannot be restricted to any interval as in the case of $X+i$.\\

      \begin{figure}[htb]
     \centering \includegraphics [width=.85\textwidth]{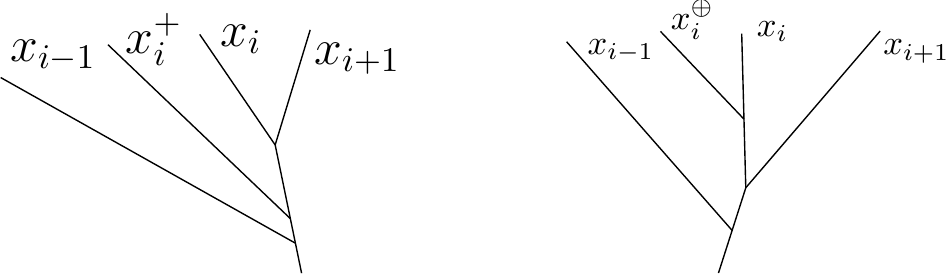}
        \end{figure}

 Notice that $X+i$ could be empty for some $X\in [\mathcal{T}]^n_{\preceq}$ and $i\in n$. For example it can be the case that $f(i)=1$ and $(x_i\cap x_{i-1})^\frown 1=x_i\cap x_{i+1}$. Observe that in this case $X+i=\emptyset$.\\
     Now we return to the inductive step. Suppose our theorem holds for all $k<n$ and consider a continuous map $c:[\mathcal{T}]^n_{\preceq}\to M$. This induces for every $i\in n$ continuous maps $c^{+}_i:[\mathcal{T}]^{n+1}_{\preceq^{+}_i}\to \{0,1\}$ defined by $c^{+}(A)=0$ if $c(A \setminus \{\alpha_i\})=c(A\setminus \{\alpha_{i+1}\})$ and $1$ if otherwise. We remind the reader here that for every $i\in n$, elements of $X+i$ are of unique order type $\{(1, \dots ,n), \preceq^+_i\}$. As a result for every $X\in [\mathcal{T}]^{n}_{\preceq}$, $X+i\subset  [\mathcal{T}]^{n+1}_{\preceq^{+}_i}$.\\
     
     By Theorem $2$ there exists a skew perfect subset $\mathcal{P}\subseteq \mathcal{T}$ and for each $i$ a constant $f_i\in \{0,1\}$ such that $c^{+}_i \upharpoonright [\mathcal{P}]^{n+1}_{\preceq^{+}_i}=\{f_i\}$. We distinguish two cases:\\
     
     $\bf{Case\, 1}:$ There exists $i\in n$ with $f_i=0$. For all $X\in [\mathcal{P}]^n_{\preceq}$ the resulting set $X:(\{0, \dots, n-1\}\setminus \{i\})$ has order type $(\{0,\dots, n-2\}, \preceq^{\star})$. Then $c$ induces naturally a continuous map $c^{+}:[\mathcal{P}]^{n-1}_{\preceq^{\star}}\to M$ defined by $$c^{+}(X:(\{0, \dots, n-1\} \setminus \{i\}))=c(X)$$ To make this map well defined and continuous, we may restrict to a perfect subset of $\mathcal{P}$.

  The inductive hypothesis now applies to gives us a skew perfect subset $\mathcal{P}_0\subseteq \mathcal{P}$, $I^{\star}$ and $J^{\star}$ such that for all $Z,W\in [\mathcal{P}_0]^{n-1}_{\preceq^{\star}}$ it holds: \\
     $c^{+}(Z)=c^{+}(W)$ iff $Z:I^{\star}=W:I^{\star}$ and $\widehat{D(Z)}:J^{\star}=\widehat{D(W)}:J^{\star}$.\\ Let $I\subseteq \{0,\dots, n-1\}$ and $J\subseteq C^n_{\preceq}$ be such that for any $X\in [\mathcal{P}_0]^n_{\preceq}$ one has:

     $$ (X:(\{ 0,\dots, n-1\} \setminus \{i\})):I^{\star}=X:I$$
     $$ \widehat{D(X:(\{ 0,\dots, n-1\} \setminus \{i\})}):J^{\star}=\widehat{D(X)}:J$$
     It follows that for all $X,Y\in [\mathcal{P}_0]^n_{\preceq}$ it holds:\\

         $c(X)=c(Y)$ iff $X:I=Y:I$ and $\widehat{D(X)}:J=\widehat{D(Y)}:J$.\\
         
         $\bf{Case\, 2:}$ Suppose that for all $i\in n$ it holds that $ f_i=1$. As we observed above elements of $X\oplus i$ are not of the same order type. The different order types can be ordered as follows. Condition $2$ of the definition of $X\oplus i$ implies that $x^\oplus_i\cap x_i$ can have length bigger than $x_{j-1}\cap x_j$, where $j\in n$ is such that $|x_i\cap x_{i+f(i)}|< |x_{j-1}\cap x_j|$ i.e. $|x_{j-1}\cap x_j|<|x^\oplus_i\cap x_i|$.\\

         In the case that $f(i)=-1$, set $B_i =\{ j\in n\setminus \{0\}: i\preceq j \}$. In the case that $f(i)=1$, set $B_i= \{ j\in n\setminus \{0\}: i+1\preceq j\}$. Equivalently $$B_i=\{\, j\in n: |x_{i}\cap x_{i+f(i)}| \leq |x_{j-1}\cap x_{j}|\, \}.$$ 
         Let $|B_i|=m_i$ and $B_i=(j_p)_{p\in m_i}$ be an enumeration such that if $q,p\in m_i$ with $q<p$ then $j_q\preceq j_p$, in other words $|x_{i}\cap x_{i+f(i)}| \leq |x_{{j_q}-1}\cap x_{j_q}| < |x_{{j_p}-1}\cap x_{j_p}|$. Observe that $j_0=i$ in the case of $f(i)=-1$ and $j_0=i+1$ in the case of $f(i)=1$.

         Then $0\in m_i$, corresponds to the total order $(\{1,\dots, n\}, \preceq^\oplus_{i,0})$ that satisfies $|x_{i+f(i)}\cap x_i|<|x^\oplus_i\cap x_i|< |x_{j_1-1}\cap x_{j_1}|$, or equivalently: $i+1\preceq^\oplus_{i,0} j_1$, when $x_{j_1}<_{lex} x_i$, and $ i+1\preceq^\oplus_{i,0} j_1+1$, when $x_i<_{lex} x_{j_1}$. A $k\in m_i$ corresponds to the total order $(\{1,\dots, n\}, \preceq^\oplus_{i,k})$ so that $$ |x_{{j_k}-1}\cap x_{j_k}|<|x^\oplus_{i}\cap x_i|<|x_{{j_{k+1}-1}}\cap x_{j_{k+1}}|$$
         We remind the reader here that the $i\in ( \{1, \dots, n-1\},\preceq)$ becomes $i+1$ in $( \{1, \dots, n\},\preceq^\oplus_{i,k})$, for every $k\in m_i$. The $i \in (\{1, \dots, n\},\preceq^\oplus_{i,k})$ comes from $x_i^\oplus\cap x_{i-1}$ and the $i+1\in (\{1, \dots, n\},\preceq^\oplus_{i,k})$ comes from $x_i^\oplus\cap x_{i}$.
         
         Having dealt with the different order types, we notice that the map $c$ gives rise to maps $c^{\oplus}_{i,k}:[\mathcal{P}]^{n+1}_{\preceq^{\oplus}_{i,k}}\to \{0,1\}$,defined by $c^{\oplus}_{i,k}(A)=0$ if $c(A-\{\alpha_i\})=c(A-\{\alpha_{i+1}\})$ and $1$ otherwise, where $k\in m_i=|B_i|$, $i\in n$.\\

         For every $i$ and $k\in m_i$, by Theorem $2$ there exists a perfect skew subset $\mathcal{P}_0$ and a constant $d_{i,k}\in \{0,1\}$ such that $c^{\oplus}_{i,k}\upharpoonright [\mathcal{P}_0]^{n+1}_{\leq^{\oplus}_{i,k}}= \{d_{i,k}\}$. In other words Theorem $2$, for every total order $(\{1,\dots, n\}, \preceq^{\oplus}_{i,k})$ gives us a constant $d_{i,k}$.\\

         Let $\{X_0,X_1\}\in [[\mathcal{P}_0]^n_{\preceq}]^2$, their union has cardinality $|X_0\cup X_1|=p$, order type $(\{1,\dots, p-1\},\preceq_p)$. We order $X_0\cup X_1$ lexicographically. There are subsets $I_0, I_1\subseteq \{0, \dots, p-1\}$ such that $(X_0\cup X_1):I_0=X_0$ and $(X_0\cup X_1):I_1=X_1$. Consider now an injective enumeration of all occurring $(\{1,\dots, p-1\},\preceq_j, I^j_0,I^j_1)_{j<q}$. For any $j<q$ let $c_j:[\mathcal{P}_0]^p_{\preceq_j}\to2$ be a map defined by $c_j(X)=0$ if $c(X:I^j_0)=c(X:I^j_1)$ and $1$ otherwise. 

        Notice that each of these maps is a continuous map. By a successive application of Theorem $2$ we get a skew perfect subset $\mathcal{P'}_1\subseteq \mathcal{P}_0$ with the property that $c_j$ restricted on $[\mathcal{P'}_1]^p_{\preceq_j}$ is constant. In other words, we get a skew perfect subset with the property that if $c(X)=c(Y)$ then $c(X')=c(Y')$ for all $X',Y'$ such that $X\cup Y$ and $X'\cup Y'$ are of the same order type $(\{1,\dots, p-1\},\preceq_j, I^j_0,I^j_1)$ for some $j\in q$. Repeat the above step for all possible $p$ and $q$ to get a skew perfect subset $\mathcal{P}_1$ such that $c(X)=c(Y)$ implies that $c(X')=c(Y')$ for all $X',Y'$ such that $X\cup Y$ and $X'\cup Y'$ are of the same order type and relative position. \\

         Let $I=\{i\in n: d_{i,k}=1$ for all $k\in m_i \, \}$. To define $J\subseteq C^n_{\preceq}$ consider the following:\\ For $i\notin I$ there exists an $l\in m_i$ such that $d_{i,l}=0$ and for all $g<l$ it holds that $ d_{i,g}=1$. We observe that if for $l$ we have $d_{i,l}=0$, then for all $h>l$ we also have $d_{i,h}=0$.  To see that, suppose there exist $X_0,X_1\in [\mathcal{P}_1]^{n}_{\preceq}$, that agree on all but their $i$-th element, namely $x^0_i\neq x^1_i$ and $x^0_j= x^1_j$ for $j<n,j\neq i$. Assume also that $c(X_0)\neq c(X_1)$ and $x^0_i$ is such that $X_1\cup \{x^0_i\}\in [\mathcal{P}_1]^{n+1}_{\preceq^\oplus _{i,h}}$.
         Then $(X_0\cup X_1)\in [\mathcal{P}_1]^{n+1}_{\preceq^{\oplus}_{i,h}}$  with respect the total order $(\{1,\dots, n\}, \preceq^\oplus_{i,h})$. Pick $x\in \mathcal{P}_1$, with $x<_{lex} x^0_i$ and $x<_{lex} x^1_i$, such that $X= (X_0-\{x^0_i\})\cup \{x\}, X_0, X_1\in [\mathcal{P}_1]^n_\preceq$ and $X_0\cup \{x\}, X_1\cup \{x\} \in [2^\omega]^{n+1}_{\preceq^\oplus_{i,l}}$.
         
         The fact that $d_{i,l}=0$ implies that $c(X_0)=c(X)=c(X_1)$, a contradiction. If no such $x\in \mathcal{P}_1$ can be found, we can always choose $X'_0,X'_1$ satisfying the above conditions with $(X_0\cup X_1), (X'_0\cup X'_1)$ being of the same order type and relative position, so that $x\in \mathcal{P}_1$ can be found. We also have $c(X')=c(Y')$.\\
         
          For $i\notin I$ and $l$ as above, namely the very first natural number with the property $d_{i,l}=0$, set $$J_i=\{\, x_i\upharpoonright |x_{j_{l}-1} \cap x_{j_{l}}|\, \}$$
          Notice that $ x_{j_{l}-1}\cap x_{j_{l}}$ is such that for any $x$, if $X\cup \{x\}\in X\oplus i$ and $x\upharpoonright |x_{j_{l}-1} \cap x_{j_{l}}|=x_i\upharpoonright|x_{j_{l}-1} \cap x_{j_{l}}|$, then $X\cup \{x\}\in [\mathcal{P}_1]^{n+1}_{\preceq^\oplus_{i,h}}$ for some $h\geq l$.           
          
          The fact that $X\cup \{x\}\in [\mathcal{P}_1]^{n+1}_{\preceq^\oplus_{i,h}}$, $h \geq l$, implies that $$c(X)=c((X\setminus \{x_i\})\cup \{x\}))$$ cause $d_{i,h}=0$ as well.
          
          Set $J=\bigcup_{i\notin I} J_i$.\\

         We claim that for all $X,Y\in [\mathcal{P}_1]^n_\preceq$ it holds that:\\

        $ c(X)=c(Y)$ if and only if $X:I=Y:I$ and $\widehat{D(X)}:J=\widehat{D(Y)}: J$.\\

        At first we show the implication from right to left. \\
        Let $X,Y$ be such that $X:I=Y:I$ and $\widehat{D(X)}:J=\widehat{D(Y)}:J$.  If $I=\{ 0,\dots ,n-1\}$ there is nothing to prove. Let $i\in \{0, \dots, n-1\}\setminus I$. There exists $l\in m_i$ such that $d_{i,l}=0$ and $l$ is the smallest possible integer with that property. By our assuption that $\widehat{D(X)}:J=\widehat{D(Y)}:J$, we conclude that $x_i\upharpoonright |x_{{j_i}-1} \cap x_{j_l}|= y_i\upharpoonright |y_{{j_i}-1} \cap y_{j_l}|$. Then $c(x_0,\dots, x_i,\dots, x_{n-1})=c(x_0, \dots, y_i, \dots, x_{n-1})$. To see this observe that $X\cup \{y_i\}\in [\mathcal{P}_1]^{n+1}_{\preceq^\oplus_{i,h}}$ for some $h\geq l$ and $d_{i,h}=0$ as well. The fact that $d_{i,h}=0$ implies the following: $c(X)=c((X\setminus \{x_i\})\cup\{y_i\})$. By doing this iteration for every $i \in n\setminus I$ we conclude that: $c(X)=c(Y)$.\\

        We show now the implication from left to right.\\
        Suppose that $c(X)=c(Y)$, and assume at first that $X:I\neq Y:I$. Let $i\in I$ be such that $x_i\neq y_i$ and suppose that $x_i<_{lex} y_i$. Choose now $y \in \mathcal{P}_1$ such that:
        
        \begin{enumerate}
            \item{} $x_i <_{lex} y <_{lex}y_i$ and $X\cup Y$ and $((X\cup Y)\setminus \{y_i\}) \cup\{ y\}$ are of the same type $( \{1,\dots, p-1\}, \preceq_j, I^j_0, I^j_1\,)$ for some $j<q$.
        
        \item{} $Y\cup \{y\}\in [\mathcal{P}_1]^{n+1}_{\preceq^\oplus_{i,k}}$ for some $k\in m_i$ .         \end{enumerate}

       These conditions imply $c(Y)=c(X)=c((Y \setminus \{y_i\})\cup\{ y\})$ which contradicts that $d_{i,k}=1$ for all $k$.

    Next consider the case that 
    $$c(X)=c(Y) \text{ and }\widehat{D(X)}:J\neq \widehat{D(Y)}:J.$$ 
    
    We claim the following.
    
    \begin{claim} 
    In our context, namely where $f_i=1$ for all $i\in n$, if $c(X)=c(Y)$, then $D(X)=D(Y)$. 
    \end{claim}
    
    \begin{proof}
       Suppose $D(X)\neq D(Y)$ and let $w\in 2^b$, $b\in \omega$, be such that $$ w\in \Delta= \{\, (D(X)\setminus D(Y))\cup (D(Y) \setminus D(X))\, \}.$$
    Assume that $w\in D(X)$ i.e. $w=x_{d-1}\cap x_{d}$, for some $d\in n$. There exists $t \in \{0,1\}$ such that $w^\frown t$ is not an initial segment of a $y_i\in Y$ for all $i\in n$. If no such a $t$ exists, then $w\in D(Y)$ as well, contradicting that $w\in \Delta$.

       Let $w^\frown t$ be an initial segment of $x_d$. Identical argument applies in the case that $w^\frown t$ is an initial segment of $x_{d-1}$.  Choose $x \in \mathcal{P}_1$ such that $(X\setminus \{x_d\})\cup \{x \}\in X+ d$ and $X\cup Y$, $((X\cup Y)\setminus \{x_d\}) \cup \{ x \}$ have the same order type, relative position. But then we have: $c(X)=c(Y)=c((X\setminus \{x_d\})\cup \{x\})$ contradicting that $f_i=1$ for all $i\in n$.\\ 
     Once more we remind the reader that if we cannot pick an $x\in \mathcal{P}_1$ directly, we can always consider $X'$, $Y'$ so that $X\cup Y$ and $X'\cup Y'$ are of the same order type and relative position, that allow us to pick such an $x$. We also have $c(X')=c(Y')$ and the above argument holds.
     \end{proof}
    
  We suppose that $c(X)=c(Y)$, $D(X)=D(Y)$ and $\widehat{D(X)}:J\neq \widehat{D(Y)}:J$ and we derive a contradiction.
  
    Our assumption $\widehat{D(X)}:J\neq \widehat{D(Y)}:J$ implies that there is an  $i\notin I $, $ x_i, y_i$ and $ k\in \omega$ such that $x_i\upharpoonright k\neq y_i\upharpoonright k$, where $x_i \upharpoonright k\in \widehat{D(X)}:J$, $y_i\upharpoonright k\in \widehat{D(Y)}:J$. Observe that the fact $D(X)=D(Y)$ implies $x_i\upharpoonright k\in \widehat{D(X)}:J$ and also $y_i\upharpoonright k\in \widehat{D(Y)}:J$ for the same $k$.
    Suppose $y_i\upharpoonright k<_{lex}x_i \upharpoonright k$. Then pick once more $x\in \mathcal{P}_1$ so that $X\cup \{x\}\in [\mathcal{P}_1]^{n+1}_{\preceq{\oplus }_{i,g}}$ for some total order $(\{1,\dots, n\},\preceq^\oplus_{i,g})$, $g<l$, where $d_{i,l}=0$ and $l\neq 0$ is the very first natural number with that property. We require also that $X\cup Y$ and $(X\cup Y\setminus \{x_i\})\cup \{ x\}$ are of the same type $(\{1,\dots, p-1\},\preceq_j, I^j_0,I^j_1)$ for $p=| X\cup Y|, j\in q$. But then $c(X)=c(Y)=c((X\setminus \{x_i\})\cup \{x\})$, contradicting that $d_{i,g}=1$ for all $g<l$. Notice that if $l=0$, then the fact that $D(X)=D(Y)$ causes the contradiction. To see this observe that $J_i=\{ \, x_i\upharpoonright |x_i\cap x_{i+f(i)}|\, \}=\{ \, x_i\cap x_{i+f(i)}\, \} \subset D(X)$.

        Once more we remind the reader that if we cannot pick an $x\in \mathcal{P}_1$ directly, we can always consider $X'$, $Y'$ so that $X\cup Y$ and $X'\cup Y'$ are of the same order type and relative position, that allows us to pick such an $x$. Then we also have that $c(X')=c(Y')$ since $\mathcal{P}_1$ has that property by construction.\\

 \end{proof}
Therefore there is a finite list of all possible patterns that correspond to each pair $(I,J)$.

 \section{Borel equivalence relations with countable quotients}

 In the case that we consider the metric space in Theorem $4$ countable then we get the following version of our Main Theorem:

 \begin{theorem} Let $n\in \omega$ and $(\{ 1,\dots, n-1\}, \preceq)$ be a total order. Further let $\mathcal{T}\subseteq 2^{\omega}$ be a perfect subset and $M$ a countable metric space. Then for any Baire measurable mapping $c: [\mathcal{T}]^{n}_{\leq}\to M$ there exists a perfect skew subset $\mathcal{P}\subseteq \mathcal{T}$, $J\subseteq C^n_\preceq$ such that  for every $X,Y \in [\mathcal{P}]^n_{\preceq}$:

 $$c(X)=c(Y) \iff  X:J=Y:J$$
 \end{theorem}

 Therefore for a countable range the above theorem states that partitions of $[\mathcal{T}]^n_\preceq $ depend only on initial segments. Theorem $5$ extends a result of Taylor \cite{Tay} for the case of $n=2$ and $M=\omega$ with the discrete topology on it:

 \begin{theorem}(Taylor)\cite{Tay} For every Baire measurable mapping $c:[2^\omega]^2\to \omega$, there exists a perfect subset $\mathcal{T}$ of $2^\omega$ such that one of the following two statements hold:\\
 $(1)$ $c\upharpoonright [\mathcal{T}]^2$ is constant\\
 $(2)$ $c(x_0,x_1)=c(y_0,y_1)$ if and only if $$d(x_0,x_1)=d(y_0,y_1)\text{, for all }(x_0,x_1)_{<_{lex}}, (y_0,y_1)_{<_{lex}} \in [\mathcal{T}]^2$$
 \end{theorem}

 And also extends the following result of Lefmann

  \begin{theorem}(Lefmann)\cite{Le} Let $M$ be a metric space and $c:[2^\omega]^2\to M$ be a Baire-measurable mapping. Then there exists a perfect subset $\mathcal{P}\subseteq 2^\omega$ and subsets $I\subseteq \{0,1\}$ and $J\subseteq \{1\}$, with $J=\emptyset$ if $I=\{0,1\}$, such that for all $(x_0,x_i)_{<_{lex}}$, $(y_0,y_1)_{<_{lex}}\in [2^\omega]^2$ it holds $$c((x_0,x_i)_{<_{lex}})=c((y_0,y_i)_{<_{lex}})$$
  $$\text{ iff }\{x_i:i\in I\}=\{y_i:i\in I\}\text{ and }\{d(x_{j-1},x_j):j\in J\}=\{d(y_{j-1},y_j):j\in J\}$$
\end{theorem}

In his approach Lefmann did not take in the account the restrictions on the intersections, namely our set $\widehat{D(X)}\setminus D(X)$. Therefore his theorem for $n=2$ has the correct form since $\widehat{D(X)}=\emptyset$.

After the completion of this work, we came across a preprint \cite{Vu} which treats the same classification problem. We were not able to verify the approach of \cite{Vu} but we are sure that our two approaches are quite different.

\end{document}